\documentclass{amsart}

\usepackage[mathscr]{euscript}
\usepackage{epsfig}
\usepackage{latexsym}

\usepackage[all]{xy}

\usepackage{amssymb}
\usepackage{amsthm}
\usepackage{amsmath}
\usepackage{amsxtra}

.tfm
.tfm

\theoremstyle{plain}
\newtheorem{prop}{Proposition}[section]
\newtheorem{thm}[prop]{Theorem}
\newtheorem{coro}[prop]{Corollary}

\newtheorem{lemma}[prop]{Lemma}
\newtheorem{conj}{Conjecture}

\theoremstyle{definition}
\newtheorem*{defi}{Definition}

\theoremstyle{remark}

\newtheorem*{remark}{Remark}

\numberwithin{table}{section}

\DeclareMathOperator{\Hom}{Hom}

\DeclareMathOperator{\sign}{sign}

\newcommand{\Ot}{{\tilde{\Om}}}
\newcommand{\Om}{{\mathscr{O}}}
\newcommand{\Os}{\Om'}
\newcommand{\Ops}{\Om^+}
\newcommand{\Ons}{\Om^-}

\def\RR{\mathbb R}

\def\<#1>{{\left\langle{#1}\right\rangle}}

\def\Z{{\mathbb Z}}             
\def\Q{{\mathbb Q}}             
\def\R{{\mathbb R}}             

\def\set#1{{\left\{{\def\st{\;:\;}#1}\right\}}}
\def\num#1{{\#{#1}}}
\def\numset#1{{\num{\set{#1}}}}

\def\tkro#1#2{\left(\tfrac{#1}{#2}\right)}             
\let\kro\dkro

\def\tf#1#2{{^{#1}\!/\!_{#2}}}

\def\H{{B}}                     
\def\Hx{{\H^\times}}            
\def\Hpx{{\H_p^\times}}         

\def\O{{R}}           
\def\Op{{\O_p}}                 
\def\Oq{{\O_q}}                 
\def\Odual{{\O^\ast}}           
\def\Ol(#1){{\mathop{\O_l}(\id{#1})}}                 
\def\Or(#1){{\mathop{\O_r}(\id{#1})}}                 
\def\Orp(#1){{\mathop{\O_r}(\idp{#1})}}               
\def\Otern(#1){{\mathop{\O^0_r}(\id{a})}}    

\def\id#1{{\mathfrak{#1}}}      
\def\idp#1{{\id{#1}_p}}         
\def\idq#1{{\id{#1}_q}}         
\def\cls#1{{[\id{#1}]}}         
\def\clsv#1{{[\id{#1}]^\vee}}   
\def\Cls#1{{\left[{#1}\right]}} 

\def\vec#1{{\mathbf{#1}}}       
\def\map#1{{{\mathit{#1}}}}              
\def\quat#1{{{\mathit{#1}}}}             
\def\quatp#1{\quat{#1}_p}     
\def\quatpconj#1{{\conjugate{\quatp{#1}}}}  

\def\conjugate#1{{\overline{#1}}}
\DeclareMathOperator{\norm}{{\mathscr N}}
\DeclareMathOperator{\trace}{{\mathrm{Tr}}}
\DeclareMathOperator{\normx}{{\Delta}}

\DeclareMathOperator{\disc}{D}
\DeclareMathOperator{\level}{N}
\DeclareMathOperator{\twolevel}{L}

\DeclareMathOperator{\omegax}{\Omega}
\def\charx#1{{\chi_{#1}}}

\DeclareMathOperator{\I}{{\mathscr I}}
\DeclareMathOperator{\Ix}{{\widetilde{\I}}}
\DeclareMathOperator{\Iv}{{\I^\vee}}
\DeclareMathOperator{\M}{{\mathscr M}}

\DeclareMathOperator{\Mv}{{\M^\vee}}
\let\MR\undefined
\DeclareMathOperator{\MR}{{\M_{\RR}}}

\def\T_#1(#2){{\mathop{\mathscr T}\nolimits_{#1}(\id{#2})}}
\def\TO(#1)_#2(#3){{\mathop{\mathscr T}\nolimits^{#1}_{#2}(\id{#3})}}
\def\t{{\mathop{t}\nolimits}}
\def\Hecke{{\phi}}

\def\HeckeRing{\mathbb{T}}

\def\A_#1(#2){{\mathop{\mathscr A}\nolimits_{#1}(\id{#2})}}
\def\Ax_#1(#2){{\mathop{\widetilde{\mathscr A}}\nolimits_{#1}(\id{#2})}}
\def\a{{\mathop{\map{a}}\nolimits}}
\def\Gross{{\Theta}}

\def\e#1{\vec{e}_{{#1}}}

\def\localconstant#1(#2){\varepsilon_{#1}({#2})}
\def\sign#1(#2){\epsilon_{#1}({#2})}

\def\charp{{\varkappa_p}}

\begin{document}

\title[Shimura correspondence for level $p^2$]%
      {Shimura correspondence for level $p^2$
       and the central values of L-series}
\author{Ariel Pacetti}
\address{Departamento de Matem\'atica, Universidad de Buenos Aires,
         Pabell\'on I, Ciudad Universitaria. C.P:1428, Buenos Aires, Argentina}
\email{apacetti@dm.uba.ar}
\thanks{The first author was supported by a CONICET grant}
\author{Gonzalo Tornar\'{\i}a}
\address{Department of Mathematics, University of Texas at Austin, Austin, TX 78712}
\email{tornaria@math.utexas.edu}
\thanks{The authors would like to thank the ``Institut Henri
      Poincar\'e'' and the organizers of the trimestre ``M\'ethodes
      explicites en th\'eorie des nombres,''
      where part of this work was written}
\begin{abstract} 
Given a weight $2$ and level $p^2$ modular form $f$, we construct
two weight $\tf{3}{2}$ modular forms (possibly zero) of level
$4p^2$ and non trivial character mapping to $f$ via the Shimura
correspondence. Then we relate the coefficients of the constructed forms to the
 central value of the L-series of certain imaginary quadratic twists of
$f$. Furthermore, we give a general framework for our construction that
applies to any order in definite quaternion algebras, with which one
can, in principle, construct weight $\tf{3}{2}$ modular forms of any level,
provided one knows how to compute ideal classes representatives.
\end{abstract}
\maketitle

\section*{introduction}
The theory of modular forms of half integral weight was developed by
 Shimura in \cite{Shimura}. There he defined a map known as the ``Shimura
 correspondence'' that associates to a half integral weight modular form
 (eigenvalue of the Hecke operators) an integral weight modular form,
 and raised the question of the meaning of their Fourier
 coefficients. Later Waldspurger related these Fourier coefficients to
 the central values of twisted L-series for the integral weight
 modular form (see \cite{Waldspurger}). In \cite{Gross} Gross gave
 (under some restrictions) an explicit method to construct, given a
 weight $2$ and level $p$ modular form $f$, a weight $\tf{3}{2}$ modular
 form (of level $4p$ and trivial character) mapping to $f$ via the
 Shimura correspondence. \par

In the first section of this paper we revise the ideal theory of
 quaternion algebras and Hecke operators acting on them.  We then
 generalize the correspondence of Gross to any order in a definite
 quaternion algebra, obtaining a Hecke linear correspondence from the
 quaternary theta series associated to ideals to the ternary theta
 series associated to their right orders.\par

In the second section we construct certain orders ``of level $p^2$''
 for which the correspondence yields modular forms of level $4p^2$ and
 show how to construct ideal classes representatives from
 representatives for the maximal order. In this way we improve the
 speed of the algorithm; moreover, the matrices that we need to
 diagonalize are much smaller. For instance, we have computed some of the weight
 $\tf{3}{2}$ modular forms for $p$ up to $500$, i.e. corresponding to
 modular forms of weight $2$ and level $p^2$ up to $250000$ (see~\cite{ell-data}). In the
 third section we give an example of how to construct these ideals
 for the case $p=7$.\par

In the fourth section we explain the relation between the Fourier
 coefficients of weight $\tf{3}{2}$ modular forms that are obtained using
 the methods of \S1 and \S2, and the central values of the L-series of
 the corresponding modular form $f$ of weight $2$. We conjecture a
 precise formula for this relation, similar to Gross' formula for the
 level $p$ case. Furthermore, we show how our conjecture implies an
 easy criteria to decide when the constructed modular forms are zero:
 essentially when $L(f,1)=0$.\par

Examples of our method, as well as an application to computing the
central values for \emph{real} quadratic twists, were presented at the
workshop ``Special Week on Ranks of Elliptic Curves and Random Matrix
Theory'' held at the Isaac Newton Institute, and can be found at
\cite{Pa-To}.


\section{Orders in quaternion algebras and Shimura correspondence}
\label{sec:quatshimura}

Let $\H$  be a definite quaternion algebra over $\Q$.
For $\quat{x}\in\H$ we denote $\norm\quat{x}$ and $\trace\quat{x}$
the \emph{reduced norm} and \emph{reduced trace} of $\quat{x}$,
respectively.
The \emph{norm} of a lattice $\id{a}$ is defined as
$\norm\id{a}:=\gcd\set{\norm\quat{x}\st\quat{x}\in\id{a}}$.
We equip $\id{a}$ with the quadratic form
$\norm_{\id{a}}(\quat{x}):=\norm\quat{x}/\norm\id{a}$,
which is primitive; its determinant is a square,
and we denote its positive square root by $\disc(\id{a})$.
In particular, when $\O\subseteq\H$ is an order, $\disc(\O)$
is its \emph{reduced discriminant}.
The subscript $p$ will denote localization at $p$, namely
$\idp{a} := \id{a} \otimes \Z_p$.

If $\O$ is an order in $\H$, we let $\Ix(\O)$ be the set of \emph{left
$\O$-ideals}, i.e. the set of lattices $\id{a}\subseteq\H$ such that
$\idp{a}=\Op\quatp{x}$ for every prime $p$, with $\quatp{x}\in\Hpx$.
Two left ideals $\id{a},\id{b}\in\Ix(\O)$ are in the same class if
$\id{a}=\id{b}\quat{x}$, with $\quat{x}\in\Hx$; we write $\cls{a}$ for
the class of $\id{a}$. The set of all left $\O$-ideal classes, which
we denote by $\I(\O)$, is known to be finite.

\subsection{The height pairing}
Let $\M(\O)$ be the free $\Z$-module with basis $\I(\O)$.
We define \emph{the height paring} by
\[
  \<\cls{a},\cls{b}> := \tfrac{1}{2}\numset{\quat{x}\in\Hx \st \id{a}\quat{x} = \id{b}} =
  \begin{cases}
     \tfrac{1}{2}\num{\Or(a)^\times} &  \text{ if $\cls{a}=\cls{b}$,} \\
     0                   &  \text{ otherwise,}
  \end{cases}
\]
where $\Or(a)$ is the right order of $\id{a}$, namely
\[
   \Or(a) := \set{\quat{x}\in\H \st \id{a}\quat{x}\subseteq\id{a}}.
\]
The height pairing induces an inner product on $\MR(\O):=\M(\O)\otimes\R$;
note that $\I(\O)$ is an orthogonal basis of this space.

The \emph{dual lattice}
$\Mv(\O):=\set{\vec{v}\in\MR(\O) \st \<\vec{v},\M(\O)>\subseteq\Z}$
is spanned by the \emph{dual basis}
\[
  \Iv(\O) := \set{\clsv{a} \st \cls{a}\in\I(\O)},
\]
where $\clsv{a}:=\frac{1}{\<\id{a},\id{a}>}\cls{a}$.
We will identify $\Mv(\O)$ with $\Hom(\M(\O),\Z)$;
indeed, a vector $\vec{v}\in\Mv(\O)$ defines a map
$\<\vec{v},\cdot>:\M(\O)\rightarrow\Z$, and conversely,
a map $\map{f}:\M(\O)\rightarrow\Z$ determines a vector
\[
  \sum_{\cls{a}\in\I(\O)} \map{f}(\cls{a})\clsv{a}\in\Mv(\O).
\]
For instance, the map $\deg:\M(\O)\rightarrow\Z$ defined by
$\deg(\cls{a}):=1$ for $\id{a}\in\Ix(\O)$, corresponds to the vector
\[
  \e{0}:=\sum_{\cls{a}\in\I(\O)} \clsv{a}\in\Mv(\O).
\]



\subsection{Hecke operators}

Let $\id{a}\in\Ix(\O)$, and $m\geq 1$ an integer. We set
\[
   \T_m(a) := \set{\id{b}\in\Ix(\O) \st
        \id{b}\subseteq\id{a}, \quad
        \norm\id{b}=m\,\norm\id{a} }.
\]
The \emph{Hecke operators} $\t_m:\M(\O)\rightarrow\M(\O)$ are then defined by
\[
  \t_m\cls{a}  := \sum_{\id{b}\in\T_m(a)} \cls{b}
\]
for $m\geq 1$ and $\cls{a}\in\I(\O)$.
In addition, we define $\t_0:\M(\O)\rightarrow\M(\O)$ by
$\t_0\cls{a} := \frac{1}{2}\e{0}$.


%

\begin{lemma}\label{lem:heckenorm}
\[
   \T_m(a) = \set{\id{b}\in\Ix(\O) \st
        \id{b}\subseteq\id{a}, \quad
        [\id{a}:\id{b}] = m^2 }.
\]
Moreover, $m\id{a}\subseteq\id{b}$ for every $\id{b}\in\T_m(a)$;
in particular, $\id{b}\in\T_m(a)$ if and only if $m\id{a}\in\T_m(b)$.
\end{lemma}
\begin{proof}
Let $\id{b}\in\Ix(\O)$ such that $\id{b}\subseteq\id{a}$.
Locally, since $\idp{a}$ and $\idp{b}$ are principal, there is some
$\quatp{x}\in\Orp(a)$ such that $\idp{b}=\idp{a}\quatp{x}$.
Then
\[
   [\idp{a}:\idp{b}] =
   (\norm\quatp{x})^2 =
   (\norm\idp{b}/\norm\idp{a})^2,
\]
which proves the first statement.
For the second, note that $\quatpconj{x}\in\Orp(a)$,
and therefore
$m\idp{a}=\idp{a}\quatpconj{x}\quatp{x}
  \subseteq\idp{a}\quatp{x}=\idp{b}$.
\end{proof}

\begin{lemma}\label{lem:matrixcount}
Fix a prime $p$ and let $\quatp{x}\in M_{2\times 2}(\Z_p)$.
The set
\[
  \set{\quatp{y}\in M_{2\times 2}(\Z_p) \st \det\quatp{y}=p, \quad
       \quatp{x}\quatp{y}^{-1} \in M_{2\times 2}(\Z_p)}
\]
is invariant under left multiplication by unimodular
matrices, and the number of orbits for this action is
$0$ if $p\nmid\det\quatp{x}$,
$1+p$ if $\quatp{x}\in p M_{2\times 2}(\Z_p)$,
and $1$ otherwise.
\end{lemma}
\begin{proof}
\def\pmat#1#2#3#4{\left(\begin{smallmatrix}{#1}&{#2}\\{#3}&{#4}\end{smallmatrix}\right)}
A well known set of representatives of
$\set{\quatp{y}\in M_{2\times 2}(\Z_p) \st \det\quatp{y}=p}$ modulo left
multiplication by unimodular matrices is
\[
\set{\pmat{p}{0}{0}{1}, \pmat{1}{0}{0}{p},
     \pmat{1}{1}{0}{p}, \dotsc, \pmat{1}{p-1}{0}{p}}.
\]
An easy calculation shows that
\begin{align*}
\pmat{a}{b}{c}{d}\pmat{p}{0}{0}{1}^{-1} \in M_{2\times 2}(\Z_p)
  & \Longleftrightarrow a\equiv c\equiv 0\pmod{p}, \\
\pmat{a}{b}{c}{d}\pmat{1}{i}{0}{p}^{-1} \in M_{2\times 2}(\Z_p)
  & \Longleftrightarrow b-ai\equiv d-ci\equiv 0\pmod{p},
\end{align*}
from which the statement follows.
\end{proof}

\begin{prop}\label{prop:heckeprop}
the Hecke operators have the following properties:
\begin{enumerate}
\item $t_m$ is self-adjoint.
\item If $(m,m')=1$, then $\t_{mm'} = \t_m\t_m'$.
\item If $p\nmid\disc(\O)$ is a prime, then
      $\t_{p^{k+2}} = \t_{p^{k+1}}\t_p - p\,\t_{p^k}$.
\end{enumerate}
\end{prop}
\begin{proof}
\begin{enumerate}
\item Note that
\begin{align*}
  \<\cls{a},\t_m\cls{b}>
  & = \sum_{\id{c}\in\T_m(b)} \<\cls{a},\cls{c}> \\
  & = \sum_{\id{c}\in\T_m(b)} \tfrac{1}{2}
      \numset{ \quat{x}\in\Hx \st \id{a}\quat{x} = \id{c} } \\
  & = \tfrac{1}{2}
      \numset{ \quat{x}\in\Hx \st \id{a}\quat{x} \in \T_m(b) }.
\end{align*}
By the last part of Lemma~\ref{lem:heckenorm}, this equals
\begin{align*}
\phantom{\<\cls{a},\t_m\cls{b}>}
  & = \tfrac{1}{2}
      \numset{ \quat{x}\in\Hx \st m\id{b} \in \T_m(a\quat{x}) } \\
  & = \tfrac{1}{2}
      \numset{ \quat{x}\in\Hx \st m\id{b}\quat{x}^{-1} \in \T_m(a) },
\end{align*}
and as before the latter is $\<\t_m\cls{a},\cls{b}>$.
\item For any $\id{c}\in\T_{mm'}(a)$, there is a unique
      $\id{b}\in\T_{m'}(a)$ such that $\id{c}\in\T_{m}(b)$,
      namely $\idp{b}=\idp{a}$ for $p\nmid m'$, and
      $\idp{b}=\idp{c}$ for $p\nmid m$.
\item Let $\id{c}\in\T_{p^{k+2}}(a)$.
      Locally, $\idp{c}=\idp{a}\quatp{x}$ for some $\quatp{x}\in\Orp(a)$
      with $\norm{\quatp{x}}=p^{k+2}$.
      Since $p\nmid D$, we can identify $\Orp(a)$ with $M_{2\times 2}(\Z_p)$,
      and use Lemma~\ref{lem:matrixcount} to count the number
      of $\id{b}\in\T_p(a)$ such that $\id{c}\in\T_{p^{k+1}}(b)$.
      Indeed, any $\id{b}\in\T_p(a)$ will be given by
      $\idp{b}=\idp{a}\quatp{y}$, where $\quatp{y}\in\Orp(a)$ is
      such that $\norm\quatp{y}=p$, and the condition
      $\id{c}\in\T_{p^{k+1}}(b)$ is equivalent to
      $\quatp{x}\quatp{y}^{-1}\in\Orp(a)$.
      \par
      Thus, if $\quatp{x}\not\in p\Orp(a)$ there is a unique
      such $\id{b}$, while for $\quatp{x}\in p\Orp(a)$ there are $1+p$
      of them. But $\quatp{x}\in p\Orp(a)$ if and only if
      $\id{c}=p\id{c'}$ for some $\id{c'}\in\T_{p^k}(a)$, and the
      formula follows since $\cls{c}=\cls{c'}$.
\end{enumerate}
\end{proof}
It follows from this Proposition
that the Hecke operators $t_m$ with $(m,\disc(\O))=1$ generate a
commutative ring $\HeckeRing^0$ of self-adjoint operators;
by the spectral theorem $\M_\R(\O)$ has
an orthogonal basis of eigenvectors for $\HeckeRing^0$.

\subsection{Modular forms of weight $2$}
The following construction will show that there is a correspondence
between $\M(\O)$ and  modular forms of weight $2$ and level
\[
   \twolevel(\O):=\norm(\Odual)^{-1},
\]
where $\Odual:=\set{\quat{x}\in\H\st\trace(x\O)\subseteq\Z}$.
Indeed, we will exhibit a $\HeckeRing^0$-linear map
\[
   \Hecke : \Mv(\O) \otimes_{\HeckeRing^0} \M(\O)
      \longrightarrow M_2(\twolevel(\O)),
\]
where $M_2(\twolevel)$ is the space of modular forms of weight $2$,
level $\twolevel=\twolevel(\O)$ and trivial character,
with $\t_n$ acting on this space by the Hecke operator $T(n)$.
We remark that this map is \emph{not} in general surjective.
\begin{defi}
Let $\vec{v}\in\Mv(\O)$ and $\vec{w}\in\M(\O)$. We set
\[
  \Hecke(\vec{v},\vec{w})
     := \sum_{m\geq 0} \<\vec{v},\t_m\vec{w}> q^m
      = \frac{\deg\vec{v}\cdot\deg\vec{w}}{2}
      + \sum_{m\geq 1} \<\vec{v},\t_m\vec{w}> q^m.
\]
\end{defi}

\begin{prop}
$\Hecke(\vec{v},\vec{w})$ is a weight $2$ modular form of level
$\twolevel(\O)$ and trivial character.
Moreover,
\[
   \Hecke(\vec{v},\vec{w})_{|T(n)}
     = \Hecke(\t_n\vec{v},\vec{w})
     = \Hecke(\vec{v},\t_n\vec{w}),
\]
for any $n\geq 1$ such that $(n,\disc(\O))=1$.
In particular, for any eigenvector $\vec{v}\in\M(\O)$ for
$\HeckeRing^0$, 
the modular form $\Hecke(\vec{v},\vec{v})$
is an eigenform for $\HeckeRing^0$.
\end{prop}
\begin{proof}
Since
\[
   \<\cls{a},\t_m\cls{b}>
    = \tfrac{1}{2}\numset{ \quat{x}\in\H \st \id{a}\quat{x}\in\T_m(b) }
    = \tfrac{1}{2}\numset{ \quat{x}\in\id{a}^{-1}\id{b}
           \st \norm\quat{x}=m\norm(\id{a}^{-1}\id{b}) },
\]
and $\<\cls{a},\t_0\cls{b}>=\tfrac{1}{2}$,
it follows that
\[
   \Hecke(\cls{a},\cls{b})
    = \frac{1}{2}\sum_{\quat{x}\in\id{c}}
             q^{\norm_{\id{c}}(\quat{x})},
\]
is the theta series of the lattice $\id{c}=\id{a}^{-1}\id{b}$,
with the quadratic form $\norm_{\id{c}}$.
Its discriminant is a square; we claim that its level is $\twolevel(\O)$.
Indeed, locally $\idp{c}=\quatp{x}\Op\quatp{y}$ for some
$\quatp{x},\quatp{y}\in\Hpx$
(where $\idp{a}=\Op\quatp{x}^{-1}$, $\idp{b}=\Op\quatp{y}$,) and
$\norm_{\id{c}}(\quatp{x}\quatp{a}\quatp{y})=\quatp{u}\norm(\quatp{a})$
for $\quatp{a}\in\Op$, where
$\quatp{u}=\norm_{\id{c}}(\quatp{x}\quatp{y})$ is a $p$-adic unit.
Therefore,
the level of $\norm_{\id{c}}$ in $\id{c}$ is equal to the level of
$\norm$ in $\O$. The latter is, by definition, the smallest positive
integer $u$ such that $u\norm(\Odual)\subseteq\Z$, since the matrix of
$\norm$ on a given basis of $\O$ is the inverse of the matrix of
$\norm$ on the dual basis of $\Odual$; but this is just
$\twolevel(\O)=\norm(\Odual)^{-1}$, as claimed.
\par
For the second statement, in view of Proposition~\ref{prop:heckeprop},
it is enough to prove the identity for $n=p\nmid\disc(\O)$ a prime.
But
\begin{align*}
  \Hecke(\vec{v},\vec{w})_{|T(p)} & =
   \sum_{m\geq 0}
   \Bigl(\<\vec{v},\t_{mp}\vec{w}>+p\<\vec{v},\t_{m/p}\vec{w}>\Bigr)\, q^m,\\
\intertext{(where $\t_{m/p}=0$ if $p \nmid m$) and}
  \Hecke(\vec{v},\t_p\vec{w}) & =
   \sum_{m\geq 0} \<\vec{v},\t_m\t_p\vec{w}>\, q^m,
\end{align*}
and the result follows since Proposition~\ref{prop:heckeprop} implies
that $\t_{mp}+p\,\t_{m/p}=t_m\,t_p$.
\end{proof}

\subsection{Special points}
Let $\id{a}\in\Ix(\O)$, and $-d\leq 0$ an integer,
$-d\equiv0,1\pmod{4}$.
The \emph{special points} of discriminant $-d$ for $\id{a}$ are
\[
  \Ax_d(a) := \set{\quat{x}\in\Or(a) \st \normx\quat{x}=-d},
\]
where $\normx\quat{x}:=(\trace\quat{x})^2-4\norm\quat{x}$ is the
discriminant of the characteristic polynomial of $\quat{x}$.
These sets are stable under translations by integers.
For each $d$ the set of orbits, which will be denoted by $\A_d(a)$, is
finite, and it is in bijection with any of the sets
\[
   \Ax_d,s(a) :=
   \set{\quat{x}\in\Ax_d(a) \st \trace\quat{x} = s} =
   \set{\quat{x}\in\Or(a) \st \trace\quat{x} = s, \quad
                              \norm\quat{x} = \frac{s^2+d}{4} },
\]
where $s$ is an arbitrary integer subject to the condition
$s\equiv d\pmod{2}$.

The maps $\a_d:\M(\O)\rightarrow\Z$ are then defined by
$\a_d(\cls{a}) = \num{\A_d(a)}$.
When $-d\not\equiv 0,1\pmod{4}$ we set $\a_d(\cls{a})=0$.
As before, we identify this maps with vectors
\[
  \e{d}:=\sum_{\cls{a}\in\I(\O)} \a_d(\cls{a})\clsv{a}\in\Mv(\O).
\]
This is consistent with the previous definition of $\e{0}$.

\begin{lemma}\label{lem:matrixcount2}
Fix a prime $p$ and let $\quatp{x}\in M_{2\times 2}(\Z_p)$.
The set
\[
  \set{\quatp{y}\in M_{2\times 2}(\Z_p) \st \det\quatp{y}=p, \quad
       \quatp{y}\quatp{x}\quatp{y}^{-1} \in M_{2\times 2}(\Z_p)}
\]
is invariant under left multiplication by unimodular
matrices, and the number of orbits for this action is
$1+p$ if $\quatp{x}\in \Z_p + p M_{2\times 2}(\Z_p)$, and
$1+\kro{\normx\quatp{x}}{p}$ otherwise.
\end{lemma}
\begin{proof}
\def\pmat#1#2#3#4{\left(\begin{smallmatrix}{#1}&{#2}\\{#3}&{#4}\end{smallmatrix}\right)}
As in the proof of Lemma~\ref{lem:matrixcount},
the statement follows from
\begin{align*}
\pmat{p}{0}{0}{1} \pmat{a}{b}{c}{d}\pmat{p}{0}{0}{1}^{-1} \in M_{2\times 2}(\Z_p)
  & \Longleftrightarrow c\equiv 0\pmod{p}, \\
\pmat{1}{i}{0}{p} \pmat{a}{b}{c}{d}\pmat{1}{i}{0}{p}^{-1} \in M_{2\times 2}(\Z_p)
  & \Longleftrightarrow ci^2+(a-d)i-b\equiv 0\pmod{p},
\end{align*}
since the discriminant of the quadratic
equation above is $\normx\quatp{x}$.
\end{proof}

\begin{lemma}\label{lem:heckelineal}
Let $\quat{x}\in\Ax_d(a)$.
Then
\[
  \numset{\id{b}\in\T_p(a) \st \quat{x}\in\Ax_d(b)} =
  \begin{cases}
     1+p            & \text{ if $\quat{x}\in \Z+p\Or(a)$,} \\
     1+\kro{-d}{p}  & \text{ if $\quat{x}\not\in \Z+p\Or(a)$,}
  \end{cases}
\]
for any prime $p\nmid\disc(\O)$.
\end{lemma}
\begin{proof}
If $b\in\T_p(a)$ we have $\id{b}$ and $\id{a}$ equal outside $p$,
and $\idp{b} = \idp{a}\quatp{y}$ for some $\quatp{y}\in\Orp(a)$
with $\norm\quatp{y}=p$; two such $\quatp{y}$ give the same $\id{b}$
if and only if they are in the same orbit under left multiplication by
units of $\Orp(a)$.  Since $p\nmid\disc(\O)$, we can identify $\Orp(a)$ with
$M_{2\times 2}(\Z_p)$, and Lemma~\ref{lem:matrixcount2} proves the claim.
\end{proof}

\begin{prop}\label{prop:edhecke}
For any prime $p\nmid\disc(\O)$ we have
\[
  \t_p\,\e{d} = \e{dp^2} + \tkro{-d}{p} \e{d} + p\,\e{d/p^2}.
\]
\end{prop}
\begin{remark}
Compare this with the formula for the action of the Hecke operators
for weight $\tf{3}{2}$ in terms of Fourier
coefficients~\cite[Theorem 1.7]{Shimura}.
\end{remark}
\begin{proof}
For an arbitrary $\cls{a}\in\I(\O)$, 
the left hand side evaluated at $\cls{a}$ is
\[
   \<\t_p\,\e{d},\cls{a}>
      = \<\e{d},\t_p\,\cls{a}>
      = \sum_{b\in\T_p(a)} \<\e{d},\cls{b}>
      = \sum_{b\in\T_p(a)} \a_d(\cls{b}),
\]
which just counts the number of pairs $(\id{b},\quat{x})$ such that
$\id{b}\in\T_p(a)$ and $\quat{x}\in\A_d(b)$.
Since $p\id{a}\subseteq \id{b}\subseteq\id{a}$,
it is clear that $\quat{x}\in\Or(b)$ implies that $p\quat{x}\in\Or(a)$.
We count the number of possible
pairs in each of three disjoint cases for $\quat{x}$:
\begin{enumerate}
\item $\numset{(\id{b},\quat{x}) \st \quat{x}\in \Z+p\Or(a)} =
       (1+p)\a_{d/p^2}(\cls{a})$.\\
There are $\a_{d/p^2}(\cls{a})$ such $\quat{x}$, and the count follows from
Lemma~\ref{lem:heckelineal}.
\item $\numset{(\id{b},\quat{x}) \st \quat{x}\in\Or(a),
                               \quad \quat{x}\not\in \Z+p\Or(a)} \\
       \mbox{}\hfill=\left(1+\kro{-d}{p}\right)\left(\a_d(\cls{a})-\a_{d/p^2}(\cls{a})\right)$.\\
There are $\a_d(\cls{a})-\a_{d/p^2}(\cls{a})$ such $\quat{x}$,
and the count follows from Lemma~\ref{lem:heckelineal}.
\item $\numset{(\id{b},\quat{x}) \st \quat{x}\not\in\Or(a)} =
       \a_{dp^2}(\cls{a})-\a_d(\cls{a})$.\\
There are $\a_{dp^2}(\cls{a})-\a_d(\cls{a})$ such $\quat{x}$;
the count follows now from Lemma~\ref{lem:heckelineal}
applied to $p\quat{x}$, since
$\normx(p\quat{x})=-dp^2$ and $\kro{-dp^2}{p}=0$.
\end{enumerate}
Adding up these expressions we get
\[
  \<\t_p\,\e{d},\cls{a}> = \a_{dp^2}(\cls{a}))
       + \tkro{-d}{p}\a_d(\cls{a})
       + \left(p-\tkro{-d}{p}\right)\a_{d/p^2}(\cls{a}),
\]
and the statement now follows using the fact that
$\kro{-d}{p}\a_{d/p^2}(\cls{a})=0$.
\end{proof}

\subsection{Modular forms of weight $\tf{3}{2}$}
We will use the special points to construct modular forms of weight
$\tf{3}{2}$. Let
\[
  \omegax = \omegax(\O) := \gcd\set{\normx\quat{x} \st x\in\O}.
\]
Note that $\normx{x}\leq 0$, and $\normx{x}=0$ if and only if $x\in\Q$;
thus
\[
  \normx_{\O}(x) := -\normx{x}/\omegax(\O)
\]
defines a \emph{primitive} positive definite ternary quadratic form on
the lattice $\O/\Z$, which we denote $Q_{\O}$.
More generally, if $\id{a}\in\Ix(\O)$, it defines a positive definite
ternary quadratic form on the lattice $\Or(a)/\Z$,
which we denote $Q_{\id{a}}$.

\begin{prop}
$Q_{\id{a}}$ is in the same genus as $Q_{\O}$.
In particular, $Q_{\id{a}}$ is integral and primitive.
Conversely, any ternary quadratic form in the genus of $Q_{\O}$ will
be equivalent to $Q_{\id{a}}$ for some $\id{a}\in\Ix(\O)$.
\end{prop}
\begin{proof}
The claim is that $\Or(a)/\Z$ is locally isometric to $\O/\Z$.
Indeed, $\idp{a}=\Op\quatp{x}$ for some $\quatp{x}\in\Hpx$, and thus
$\Orp(a)=\quatp{x}^{-1}\Op\quatp{x}$, inducing an isometry between
$\Op/\Z_p$ and $\Orp(a)/\Z_p$.
\par
Conversely, let $Q$ be a quadratic form in the genus of $Q_{\O}$.
By the correspondence between ternary quadratic forms and orders in
quaternion algebras (see~\cite{Llorente}), there is an order $\O'$
in $\H$ such that $Q_{\O'} \sim Q$. Moreover, since $Q$ and $Q_{\O}$
are in the same genus, it follows that $\O'$ and $\O$ are locally
conjugate, i.e. $\O'_p = \quatp{x}^{-1}\Op\quatp{x}$.
Thus the right order of the $\O$-ideal $\id{a}$ given by
$\idp{a}=\Op\quatp{x}$ will be $\O'$, and so
$Q_{\id{a}}\sim Q$.
\end{proof}

\begin{coro} $\e{d}=0$ unless $d\equiv 0\pmod\omegax$.
\end{coro}
\begin{proof}
Note that $\a_d(\cls{a})$ equals the number of representations of $d/\omegax$
by $Q_{\id{a}}$ which, being integral, represents only integers.
Thus $\omegax\nmid d$ implies $\a_d(\cls{a})=0$ for any $\id{a}$,
so that $\e{d}=0$.
\end{proof}

\begin{prop}\label{prop:ternlevel}
The level of $Q_{\id{a}}$ is $4\twolevel(\O)/\omegax(\O)$,
and its discriminant is, up to squares, $\omegax(\O)$.
\end{prop}
\begin{proof}
It is enough to find the level and discriminant for $Q_{\O}$.
Consider the map $\rho:\O/\Z\rightarrow\O$
given by $\rho(\quat{u})=\quat{u}-\conjugate{\quat{u}}$.
Note that $\norm{\rho(u)}=\omegax\cdot\normx_{\O}\quat{u}$.
Now, if $\quat{u}\in\O$ we have
$\rho(\quat{u})=2\quat{u}-\trace{u}\in\Z+2\O$, and conversely
for $t\in\Z$ we can write,
$t+2\quat{u}=t+\trace{u}+\rho(\quat{u})\in\Z+\rho(\O/\Z)$.
Thus $\Z+\rho(\O/\Z)=\Z+2\O$, where the first sum is orthogonal (with
respect to the quadratic form $\norm$). It follows that the quadratic
form $\norm$ in the lattice $\Z+2\O$ is equivalent to the quadratic
form $1+\omegax\cdot Q_{\O}$. The assertion follows now from a
straightforward calculation, since the determinant of $\Z+2\O$ is a
square, and its level is $\twolevel(\Z+2\O)=4\twolevel(\O)$.
\end{proof}

Accordingly, we define the \emph{level} and the \emph{character} of an order $\O$
to be
\[
   \level(\O) := 4\twolevel(\O)/\omegax(\O),
   \qquad\text{and}
   \qquad \charx{\O}(n) := \kro{\omegax(\O)}{n}.
\]
We will be constructing a $\HeckeRing^0$-linear map
\[
   \Gross: \M(\O) \rightarrow
   M_{\tf{3}{2}}\left(\level(\O),\charx{\O}\right),
\]
where $M_{\tf{3}{2}}(\level(\O),\charx{\O})$ is the space of modular forms of
weight $\tf{3}{2}$, level $\level=\level(\O)$ and character
$\charx{}=\charx{\O}$, with $\t_n$ acting on this space by the Hecke operator
$T(n^2)$ (see \cite{Shimura} for the definition of modular forms of half integral weight
and the Hecke operators acting on them.)

\begin{defi}
Let $\vec{v}\in\M(\O)$. We set
\[
  \Gross(\vec{v})
     := \frac{1}{2}\sum_{d\geq 0} \<\e{d},\vec{v}> q^{d/\omegax}
      = \frac{\deg\vec{v}}{2} + \frac{1}{2}\sum_{d\geq 1} \a_d(\vec{v})
          q^{d/\omegax}
\]
\end{defi}

\begin{prop}
$\Gross(\vec{v})$ is a weight $\tf{3}{2}$ modular form of level
$\level(\O)$ and character $\charx{\O}$.
It is a cusp form if and only if $\deg\vec{v}=0$.
Moreover,
\[
   \Gross(\vec{v})_{|T(n^2)} = \Gross(\t_n\vec{v}),
\]
for any $n\geq 1$ such that $(n,\disc(\O))=1$.
\end{prop}
\begin{proof}
We have
\[
   \Gross(\id{a}) = \frac{1}{2}\sum_{\quat{x}\in{\Or(a)/\Z}}
                         q^{\normx_{\O}(\quat{x})},
\]
is the theta series of the quadratic form $Q_{\id{a}}$,
and the claim on the level and character follows from
Proposition~\ref{prop:ternlevel}.
\par
For the second claim, note that $\Gross(\vec{v})$ is a linear
combination of theta series corresponding to quadratic forms in the
same genus; thus it vanishes at all the cusps if and only if it vanishes
at the $\infty$ cusp.
\par
The last statement is exactly Proposition~\ref{prop:edhecke}.
\end{proof}

\begin{prop}
If $\O' = \Z + b\O$ for some $b\in\Z$, then
$\disc(\O')=b^3\disc(\O)$, $\twolevel(\O')=b^2\twolevel(\O)$,
and $\omegax(\O')=b^2\omegax(\O)$,
hence
\[
   \level(\O')=\level(\O), \qquad\text{and}\qquad \charx{\O'}=\charx{\O}.
\]
Moreover, $\Gross(\M(\O')) = \Gross(\M(\O))$.
\end{prop}
\begin{proof}
Indeed since $\O'/\Z = b(\O/\Z)$ it is obvious that
$\omegax(\O')=b^2\omegax(\O)$ and that indeed
$Q_{\O'}=Q_{\O}$. Everything else follows easily as in the proof of
Proposition~\ref{prop:ternlevel}.
\end{proof}

\begin{defi} An order $\O$ is called \emph{primitive} if it is not of
the form $\O = \Z+b \O'$ with $b \in \Z$, $b \neq \pm 1$ and $\O'$ an
order.
\end{defi}

\begin{coro} For the purpose of constructing modular forms of weight
$\tf{3}{2}$, it is enough to consider primitive orders.
\end{coro}

\subsection{Subideals}

Let $\O, \O'$ be orders and $\id{a}$ be a left $\O$-ideal. We define

\[
\Psi^{\O}_{\O'}(\id{a}) := \set{\id{b}\in\Ix(\O') \st \id{b}\subseteq\id{a}, \quad
  \norm\id{b}=\norm\id{a} }.
\]

This induces a map $\psi^{\O}_{\O'}: \M(\O) \mapsto \M(\O')$ by
$\psi^{\O}_{\O'}([\id{a}]) := \sum _{\id{b} \in
    \Psi^{\O}_{\O'}(\id{a})} [\id{b}]$.


\begin{prop} Let $\O$ and $\O'$ be orders in $\H$. Then 
\[
\t_m \, \psi^{\O}_{\O'} = \psi^{\O}_{\O'} \, \t_m
\]
provided $(\O)_p = (\O')_p$ for all primes $p$ dividing $m$.
\end{prop}

\begin{proof} Let $\id{a} \in \Ix(\O)$. Given $\id{b} \in
  \Psi^{\O}_{\O'}(\id{a})$ and $\id{c}\in\T_m(b)$ there is a unique
  $\id{b'} \in \T_m(a)$ such that $\id{c} \in
  \Psi^{\O}_{\O'}(\id{b'})$. Namely, if $p \nmid m$, $\idp{b}=
  \idp{c}$ and $\idp{b'} = \idp{a}$. Otherwise, $\idp{a} = \idp{b}$
  and $\idp{b'} = \idp{c}$. Similarly, given $\id{b} \in \T_m(a)$ and
  $\id{c} \in \Psi^{\O}_{\O'}(\id{b})$ there is a unique $\id{b'}$
  such that $\id{b'} \in \Psi^{\O}_{\O'}(\id{a})$ and
  $\id{c}\in\T_m(b')$.
\end{proof}

In the particular case where $\O' \subset \O$ the
hypothesis of the Proposition
is equivalent to
$\gcd(m,[\O':\O]) = 1$.\\



\section{Orders of level $p^2$}
\label{sec:ordersp2}



Fix a prime $p>2$, and let $\H$ be the quaternion algebra over $\Q$
which is ramified at $p$ and $\infty$. Let $\Om$ be a maximal order in
$\H$ and let $\Ot := \set{\quat{x} \in \Om : p \mid \normx\quat{x}}$ be the
unique order of index $p$ in $\Om$. Note $\twolevel (\Ot) = \disc(\Ot)
= p^2$, but $\level(\Ot) = 4p$ since $\omegax (\Ot) = p$.

\begin{lemma}
\label{order-condition}
\emph{
  Let $\O$ be an order, and let $\quat{x}\in \H$. A necessary and sufficient
  condition for $\O'=\Z x + \O$ to be an order
  is that $x$ be integral and $x\O\subseteq\O'$.
}
\par
\end{lemma}
\begin{proof}
Clearly if $\O'$ is an order, $\quat{x}$ has to be integral,
and since $\O' \O' = \O'$, we have $\quat{x}\O \subseteq \O'$.
The other implication follows from the
fact that $\quat{x}$ being integral implies
$\quat{x}^2=\trace(\quat{x})\quat{x} -\norm(\quat{x})\in\O'$, hence
$\O'$ is closed by multiplication.
\end{proof}

\begin{prop}
  Let $L \subset \Ot$ be a lattice such that $[\Ot:L] = p$. Then
  $L$ is an order if and only if $\Z + p\Om \subset L$.
\end{prop}
\begin{proof}

  
  Let $\Os$ be an order of index $p$ in $\Ot$. Locally, $\Os$ is already maximal
  outside $p$ hence all orders are the same. At $p$ there is a unique
  maximal order since $\H$ is ramified there. Let
  $\set{u_0=1,u_1,u_2,u_3}$ be an orthogonal basis for $\Os_p$. If $\Z
  + p\Om \not \subset \Os$ then there exists $v \in \Os$ such that
  $\Om = \Z \frac{v}{p^2} +\Os$. Since $p^4 \mid \norm(v)$ we claim
  that one of the basis elements of $\Os_p$ has norm divisible by
  $p^4$.  Note that $v/p \not \in \Os$ hence $v$ can be written as $v
  = \sum a_i u_i$ with some $a_i$ not divisible by $p$. If $p^4 \nmid
  \norm(u_i)$ for $i=0,..,3$ we would have a non zero solution of
  $\norm(v) = 0 \bmod(p^4)$ which by Hensel Lemma would lift to $w_p$,
  a non zero element in $\H_p$ with $\norm(w_p) = 0$. This cannot
  happen since the quadratic form norm is anisotropic on $\H_p$. Hence
  $p^4 \mid \norm(u_i)$ for some $i$ (say $p^4 \mid \norm(u_3)$).
  Then $\Z_p \frac{u_3}{p^2} + \Os_p= \Om_p$. $\Os_p$ being an order
  and the chosen basis being orthogonal implies $u_1u_2 = ku_3$ with $k
  \in \Z_p$, i.e. $\Om = \set{1,u_1,u_2,\frac{u_1u_2}{k}}$. But $p^2
  \nmid N(u_i)$ for $i=1,2$ and $\frac{u_1u_2}{p^2k} \in \Om$ then
  $\frac{\norm(u_1)\norm(u_2)}{p^4k^2} \in \Z_p$ therefore $k \not \in
  \Z_p$ which is a contradiction.

Conversely, if $L$ is such a lattice, let $x\in L$
such that $x\not\in\Z+p\Om$, so that $L=\Z x+(\Z+p\Om)$.
But $x(\Z+p\Om)=\Z x + px\Om\subseteq\Z x+(\Z+p\Om)$,
and Lemma~\ref{order-condition} implies that $L$ is an order.
\end{proof}

The orders as in the above Proposition will be said to be the
\emph{orders of level $p^2$}, and will be denoted by $\Os$.
Note that $\twolevel(\Os) = \disc(\Os) = p^3$ and $\level(\Os) = 4p^2$.

\begin{remark}
There are no orders $\O$ with $\twolevel(\O) = p^2$ and $\level(\O) = 4p^2$.
\end{remark}

\begin{lemma}
  Let $x \in \Os$. If $x \not \in \Z+p\Om$ then $p\parallel \normx x$
  and $\kro{\normx x/p}{p}=\pm 1$ is independent of $x$.
\end{lemma}
\begin{proof}
Let $x_0 \in \Os$
be such that $\Os = \Z x_0 + (\Z + p\Om)$. Any element $x \in \Os$ is of the form $a
x_0 + v$ where $v \in \Z+p\Om$, hence $x \not \in \Z+p\Om$ if and only
if $p \nmid a$. But $\normx x = a^2 \normx x_0 + 2a (\norm (x_0v)-
\norm(x_0 \bar v)) + \normx x \equiv a^2 \normx x_0 \pmod {p^2}$, hence the
Kronecker symbol is independent of $x$.
\end{proof}

We define the character of $\Os$ to be
\[
   \sigma(\Os) := \kro{\normx x/p}{p}
\]
where $x$ is in the conditions of the Lemma.

\begin{prop}
Two orders $\Os_1$ and $\Os_2$ of level $p^2$ 
are locally conjugate if and only if $\sigma(\Os_1)=\sigma(\Os_2)$.
\end{prop}
\begin{proof}
It is clear that $\sigma(\Os)$ is an invariant by conjugation, since
$\normx(\alpha x \alpha^{-1}) = \normx x$ for any $\alpha \in
\H_p^\times$.  \par For the converse, let $x_i\in (\Os_i)_p$ ($i=1,2$)
such that $\trace(x_i)=0$ and $(\Os_i)_p=(\Z_p+p\Om_p)+\Z_p x_i$. We
can assume that $\norm(x_1)=\norm(x_2)$, since $\normx x_i = 4 \norm
(x_i)$ implies $\norm(x_1/x_2)\in(\Z_p^\times)^2$. In this case there
is an element $\alpha \in \H_p$ that sends $x_1$ to $x_2$ via
conjugation. Clearly conjugation by $\alpha$ sends $(\Os_1)_p$ onto
$(\Os_2)_p$ since the maximal order being unique at $p$, any
conjugation has to send $\Om_p$ onto $\Om_p$, and thus $\Z_p+p\Om_p$
onto $\Z_p+p\Om_p$.
\end{proof}




Define the group ${\mathscr H} := \Z_p^\times(1+p\Om_p) \backslash
\Om_p^\times$, the group 
${\mathscr G} := \Z_p^\times(1+p\Om_p) \backslash \Ot_p^\times$ and
the group ${\mathscr S} := \Z_p^\times(1+p\Om_p) \backslash
(\Os_p)^\times$. 
\begin{lemma} The groups defined above satisfy:
\begin{enumerate}
\item $\# {\mathscr H} = p^2(p+1)$, $\# {\mathscr G} = p^2$ and $\#
  {\mathscr S} = p$.
\item $ {\mathscr S} \triangleleft {\mathscr G} \triangleleft
  {\mathscr H}$.
\end{enumerate}
\end{lemma}

\begin{proof}
Let $\pi$ be a generator of the unique order of norm $p$ in $\Om_p$ (in
particular $\pi^2 = pu$ with $u \in \Z_p^\times$),
then :

\begin{equation*}
\Om_p^\times = \set{ a_0 + a_1 \pi + a_2 p+... \st a_i \in \Om_p /
  \langle \pi \rangle \text{ and } a_0 \in (\Om_p / \langle \pi
  \rangle)^ \times}
\end{equation*}
\begin{equation*}
\Ot_p^\times = \set{ \Z_p^\times + a_1 \pi + a_2 p+... \st a_i \in \Om_p /
  \langle \pi \rangle}
\end{equation*}
and
\begin{equation*}
\Z_p^\times(1+pO_p)= \set{ \Z_p^\times + a_2 p+... \st a_i \in \Om_p /
  \langle \pi \rangle \text{ and } a_0 \in (\Om_p / \langle \pi
  \rangle)^ \times}
\end{equation*}
These equations and the fact that $\Os_p \neq \Z_p + p\Om_p$ imply the first statement.
Also it is clear that $\Ot_p ^\times \triangleleft \Om_p^\times$ and
$(\Os_p)^\times \triangleleft \Ot_p^\times$ for any order $\Os$ which
proves the second statement. 
\end{proof}

Given an order $\O$ and $\mathscr{A}\subset \Ix(\O)$, we denote
$\Cls{\mathscr{A}} := \set{ \cls{a} \st \id{a} \in \mathscr{A}}$.

\begin{thm}\label{thm:idealsp2}
Let $\Om$, $\Ot$ and $\Os$ as before. Then:
\begin{enumerate}
\item $\I(\Ot) = \bigsqcup_{\cls{a} \in \I(\Om)} \Cls{\Psi_{\Ot}^\Om(\id{a})}$
    (disjoint union).
\item The set $\Psi_{\Ot}^\Om(\id{a})$ is a principal homogeneous
    space for the cyclic group ${\mathscr G} \backslash {\mathscr H}$.
\item $\I(\Os) = \bigsqcup_{\cls{b} \in \I(\Ot)}
  \Cls{\Psi_{\Os}^\Ot(\id{b})}$. Furthermore $\# \Cls{\Psi_{\Os}^\Ot(\id{b})}
  = \# \Psi_{\Os}^\Ot(\id{b}) = p$.
\item All elements in $\Psi_{\Os}^\Ot(\id{b})$ have the same right
  order for any $\cls{b} \in \I(\Ot)$.
\item $\Psi_{\Os}^\Ot(\id{b})$ is a principal homogeneous space for
  ${\mathscr S} \backslash {\mathscr G}$.
\end{enumerate}
\end{thm}

\begin{proof} Points $(1)$ and $(2)$ are subsection 3.3 of
  \cite{Ville}, where the action of $\quatp{x}\in\Hpx$ on an ideal $\id{b}$ is given
  by right multiplication by the adele
\[
(\quat{x}_q) =
\begin{cases}
\quatp{x} & \text{ if $q=p$,}\\
1 & \text{ if $q\neq p$.}
\end{cases}
\]
 Note that ${\mathscr G} \backslash {\mathscr H}$ might
  not act freely on $\Cls{\Psi_{\Ot}^\Om(\id{a})}$.
  Indeed, one can see that the stabilizer of this action
  is a subgroup of order $\<\cls{a},\cls{a}>$.
  
  In $(3)$ the union is clearly disjoint: if $\id{c}_i \in
  \Psi_{\Os}^\Ot(\id{b}_i)$ with $\id{b}_i \in \I(\Ot)$ for $i=1,2$
  and $\id{c}_1 = \id{c}_2 \quat{x}$ for some $\quat{x} \in \H$, then
  $\id{b}_1 \quat{x} = \Ot \id{c}_1 \quat{x} = \Ot \id{c}_2 = \id{b}_2$.
  Now let $\id{c}$ be any left $\Os$-ideal.
  Then $\Ot\id{c}$ is a left $\Ot$-ideal, and
  $\norm\id{c}=\norm(\Ot\id{c})$, hence
  $\id{c}\in\Psi_{\Os}^\Ot(\Ot\id{c})$.
  
  Hence we are led to prove that $\Psi_{\Os}^\Ot(\id{b})$ has exactly
  $p$ elements, all non equivalent. Clearly if $\id{c} \in
  \Psi_{\Os}^\Ot(\id{b})$, $[\id{b}:\id{c}] = p$. 

\begin{lemma}
The set $\Psi_{\Os}^\Ot(\id{b})$ is non empty.
\end{lemma}
\begin{proof} Let $\id{b}$ be a left $\Ot$ ideal, say 
  $\idq{b} = \Ot_q\quat{x}_q$, then
  the lattice $\id{c}$ given locally by $\idq{c} = \Os_q\quat{x}_q$ is in
  $\Psi_{\Os}^\Ot(\id{b})$.
\end{proof}

Consider the left action of ${\mathscr G}$ on
$\Psi_{\Os}^\Ot(\id{b})$ defined by $g \in {\mathscr{G}}$ on an ideal $\id{c} \in
\Psi_{\Os}^\Ot(\id{b})$,  say $\id{c}_q = \Os_q \quat{x}_q$ locally, by 
\[
\bigl(g\cdot\id{c}\bigr)_q = \begin{cases}
\Os_p\, g \quat{x}_p & \text{ if $q=p$,}\\
\Os_q \quat{x}_q = \id{c}_q & \text{ if $q\neq p$.}
\end{cases}
\]
Since $(\Os_p)^\times \triangleleft \Ot_p^\times$ it is easy to check
that this action is well defined and that the stabilizer of this
action is ${\mathscr S}$. Furthermore if $\id{c}, \id{d} \in
\Psi_{\Os}^\Ot(\id{b})$, then $\Ot \id{c} = \Ot \id{d}$. Hence,
if $\id{c}_p = \Os_p \quatp{x}$ and $\id{d}_p = \Os_p \quatp{y}$,
there exists $g \in \Ot_p^\times$ such that $\quatp{x} = g \quatp{y}$,
thus $g\cdot\id{d} = \id{c}$, i.e. $\Psi_{\Os}^\Ot(\id{b})$ is a
principal homogeneous space for ${\mathscr S} \backslash {\mathscr
  G}$. In particular the number of elements in
$\Psi_{\Os}^\Ot(\id{b})$ is $p$ and they all have the same right
order. Two elements cannot be equivalent since there are no units in
$\Ot$ other than $\pm 1$ and all elements in $\Psi_{\Os}^\Ot(\id{b})$
have the same norm.
\end{proof}

A left $\Ot$-ideal $\id{b}\in\Psi_{\Ot}^\Om(\id{a})$ can be computed
using Lemma 2 in~\cite{Ville}. Acting on the right of $\id{b}$ by
representatives of ${\mathscr G} \backslash {\mathscr H}$ we obtain
all of $\Psi_{\Ot}^\Om(\id{a})$. Repeating for all $\id{a}$ in a set
of representatives of $\I(\Om)$, we can obtain a set of
representatives for $\I(\Ot)$. 

\begin{lemma}
Let $\id{m} = \id{m}_{\Os} := \set{ \quat{x} : \quat{x} \Ot \subset \Os }$.
Then $\id{m}_{\Os}$ is a bilateral $\Ot$-ideal of index $p^2$ in $\Ot$.
\end{lemma}
\begin{proof}
Since $\idq{m} = \Ot_q$ for all primes $q \not = p$ it is clear that
$\id{m}$ is bilateral. Also from the chain $\Op \supsetneq \Ot_p \supsetneq
\Os_p \supsetneq \id{m}_p \supsetneq p\O_p$ we see that $\id{m}$ has index $p^2$
in $\Ot$. 
\end{proof}

\begin{prop}\label{prop:subsubideal}
  $\id{m}_{\Os} \id{b} \subset \id{b}$ with index $p^2$ for any left
  $\Ot$-ideal $\id{b}$.  Furthermore,
  \[
    \Psi_{\Os}^\Ot(\id{b}) = \set{ \id{c} \st \id{m}_{\Ot} \id{b} \subsetneq
  \id{c} \subsetneq \id{b} \text{ and } \norm \id{c} = \norm
  \id{b}}.
  \]
\end{prop}
\begin{proof}
  The first claim follows directly from the Lemma.
  Thus, the number of lattices $\id{c}$ such that $\id{m}_{\Ot} \id{b} \subsetneq
  \id{c} \subsetneq \id{b}$ is $p+1$. The left $\Om$-ideal
  corresponding to the different of $\Om$ times $\id{b}$ is among
  these $p+1$ lattices but has norm $p \norm \id{b}$, then the set on
  the right is empty or has exactly $p$ elements. Hence we are led to
  prove that there exists a left $\Os$-ideal $\id{c} \in
  \Psi_{\Os}^\Ot(\id{b})$ such that $\id{m}_{\Ot} \id{b}
  \subsetneq \id{c} \subsetneq \id{b}$ and $\norm \id{c} = \norm
  \id{b}$. Let $\id{c}_q = \id{b}_q$ for all primes $q \neq p$ and if
$\idp{b} = \Ot_p\quatp{x}$, let $\idp{c} = \Os_p\quatp{x}$, then $\id{c}$
has the required properties.
\end{proof}

\begin{remark} We can obtain representatives for $\I(\Os)$ by applying
this Proposition to each one of a set representatives for $\I(\Ot)$.
\end{remark}




\section{Example: prime $p = 7$}
Let $\H=\H(-1,-7)$,
the quaternion algebra ramified precisely at $7$ and $\infty$. This
quaternion algebra has class number $1$ (hence type number $1$ also). 
A maximal order is given by
\begin{align*}
   \Om=\<1,i,\frac{1+j}{2},\frac{i+k}{2}>. \\
\end{align*}
Its index $p$ suborder is given by
\[
  \Ot = \<1,7i,\frac{1+j}{2},\frac{7i+k}{2}>;
\]

A generator of ${\mathscr G} \backslash {\mathscr H}$ is the element
$\frac{1+2i+j}{2}\in\Om_p^\times$.
Hence, we get all the left $\Ot$-ideals by repeatedly acting on
$\Ot$ by this element.

Let $\id{b}$ be a lattice in $\H$, and let $\quatp{x}\in\Hpx$.
We denote $\id{b}\star\quatp{x}$ the action,
by right multiplication, of the adele
\[
(\quat{x_q}) = \begin{cases}
\quatp{x} & \text{ if $p=q$,}\\ 1 & \text{ if $p \neq q$,}
\end{cases}
\]
on the lattice $\id{b}$.  Let $\O$ be the right order of $\id{b}$.
We can compute this action as follows:
\begin{enumerate}
\item Let $k$ be the smallest integer such that
      $p^k \Op\subseteq\Op\quatp{x}$ (so that
      $p^k \idp{b} \subseteq \idp{b}\quatp{x}$.)\\
      For instance, let $k=t+s$, where $t$ is the valuation at $p$ of
      $\norm\quatp{x}$, and $s$ is the smallest integer such that
      $p^s \quatp{x}\in \Op$.
\item Let $\quat{y}\in\O[p^{-1}]$ such that
      $\quat{y}-\quatp{x} \in p^k \Op$.\\
      For instance, write $\quatp{x}$ in a basis of $\O$ with
      coefficients in $\Q_p$, and then reduce the coefficients modulo
      $p^k$.
\item It now follows from a local computation that $\id{b}\star\quatp{x}=p^k\id{b} + \id{b}y$. \\
      Indeed, at $q\neq p$ we have $\quat{y}\in\Oq$, and
      $\idq{b}=p^k\idq{b}$.\\
      At $p$ it follows from (2), since
      $p^k\idp{b}\subseteq\idp{b}\quatp{x}$.
\end{enumerate}

As an example, consider the case $\id{b}_1=\Ot$
and $\quatp{x}=\frac{1+2i+j}{2}\in\Hpx$, with $p=7$.
The right order of $\id{b}_1$ is $\O=\Ot$, and we can take $k=1$.
Note that $\quatp{x}\in\O[1/7]$, hence we can use
$\quat{y}=\frac{1+2i+j}{2}\in\Hx$.
It is now easy to compute $\id{b}_2 = 7\id{b}_1 + \id{b}_1 \quat{y}$.

Repeating, we obtain $\Psi_{\Ot}^\Om(\id{\Om}) = \set{\id{b}_i : 1 \le i \le 8}$, where
the left $\Ot$-ideals $\id{b}_i$ and their characters are shown in
Table~\ref{table:psi:7}. 
The ideals $\id{b}_i$ and $\id{b}_{4+i}$ for $1 \le i \le 4$ are equivalent
(since $\<\cls{a},\cls{a}>=2$),
hence
\[
\I(\Ot) = \set{\id{b}_i : 1 \le i \le 4}.
\]
\begin{table}
\begin{tabular}{|l||l|}
\hline
$\Ot$-subideals & $\chi$ \\
\hline\hline & \\[-2.5ex]
$\id{b}_1 = \<1,7i,\frac{1+j}{2},\frac{7i+k}{2}>$ & $+$\\
$\id{b}_2 = \<7,4+i,\frac{7+j}{2},\frac{4+i+k}{2}>$ & $-$\\
$\id{b}_3 = \<7,1+i,\frac{7+j}{2},\frac{8+i+k}{2}>$ & $+$\\
$\id{b}_4 = \<7,2+i,\frac{7+j}{2},\frac{2+i+k}{2}>$ & $-$\\
$\id{b}_5 = \<7,i,\frac{7+j}{2},\frac{i+k}{2}>$ & $+$\\
$\id{b}_6 = \<7,5+i,\frac{7+j}{2},\frac{12+i+k}{2}>$ & $-$\\
$\id{b}_7 = \<7,6+i,\frac{7+j}{2},\frac{6+i+k}{2}>$ & $+$\\
$\id{b}_8 = \<7,3+i,\frac{7+j}{2},\frac{10+i+k}{2}>$ & $-$\\[.5ex]
\hline
\end{tabular}
\caption{Table of left $\Ot$-subideals.
          \label{table:psi:7}}
\end{table}

We fix two index $p$ suborders of $\Ot$
\begin{align*}
  \Ops & = \<1,7i,\frac{1+j}{2},\frac{7i+7k}{2}>, \\
  \Ons & = \<1,7i,\frac{1+7j}{2},\frac{1+7i+5j+k}{2}>
\end{align*}
in the $+$ and $-$ genus respectively.  Table~\ref{table:O:7} shows
the subideals under each $\id{b_i}$ for $\Ops$ (respectively $\Ons$),
for $i=1,\dotsc,4$.
\begin{table}
\def\<#1>{\bigl\langle#1\bigr\rangle}
\begin{tabular}{|l|l|l|}
\hline
$\Ot$-ideals & $\Ops$-subideals & $\Ons$-subideals \\
\hline\hline & & \\[-4ex]
& $\<1,7i,\frac{1+j}{2},\frac{7i+7k}{2}>$ & $\<1,7i,\frac{1+7j}{2},\frac{1+7i+5j+k}{2}>$ \\
& $\<7,7i,\frac{7+j}{2},\frac{2+7i+k}{2}>$ & $\<7,7i,\frac{1+j}{2},\frac{2+7i+k}{2}>$ \\
& $\<7,7i,\frac{7+j}{2},\frac{4+7i+k}{2}>$ & $\<7,7i,\frac{3+j}{2},\frac{6+7i+k}{2}>$ \\
$\id{b}_1$
& $\<7,7i,\frac{7+j}{2},\frac{6+7i+k}{2}>$ & $\<7,7i,\frac{5+j}{2},\frac{10+7i+k}{2}>$ \\
& $\<7,7i,\frac{7+j}{2},\frac{8+7i+k}{2}>$ & $\<7,7i,\frac{9+j}{2},\frac{4+7i+k}{2}>$ \\
& $\<7,7i,\frac{7+j}{2},\frac{10+7i+k}{2}>$ & $\<7,7i,\frac{11+j}{2},\frac{8+7i+k}{2}>$ \\
& $\<7,7i,\frac{7+j}{2},\frac{12+7i+k}{2}>$ & $\<7,7i,\frac{13+j}{2},\frac{12+7i+k}{2}>$ \\[1ex]
\hline & & \\[-4ex]
& $\<7,4+i,\frac{7+7j}{2},\frac{11+i+3j+k}{2}>$ & $\<7,4+i,\frac{7+7j}{2},\frac{4+i+k}{2}>$ \\
& $\<7,7i,\frac{1+2i+j}{2},\frac{4+i+k}{2}>$ & $\<7,7i,\frac{1+2i+j}{2},\frac{7i+k}{2}>$ \\
& $\<7,7i,\frac{3+6i+j}{2},\frac{12+3i+k}{2}>$ & $\<7,7i,\frac{3+6i+j}{2},\frac{7i+k}{2}>$ \\
$\id{b}_2$
& $\<7,7i,\frac{5+10i+j}{2},\frac{6+5i+k}{2}>$ & $\<7,7i,\frac{5+10i+j}{2},\frac{7i+k}{2}>$ \\
& $\<7,7i,\frac{9+4i+j}{2},\frac{8+9i+k}{2}>$ & $\<7,7i,\frac{9+4i+j}{2},\frac{7i+k}{2}>$ \\
& $\<7,7i,\frac{11+8i+j}{2},\frac{2+11i+k}{2}>$ & $\<7,7i,\frac{11+8i+j}{2},\frac{7i+k}{2}>$ \\
& $\<7,7i,\frac{13+12i+j}{2},\frac{10+13i+k}{2}>$ & $\<7,7i,\frac{13+12i+j}{2},\frac{7i+k}{2}>$ \\[1ex]
\hline & & \\[-4ex]
& $\<7,1+i,\frac{7+7j}{2},\frac{8+i+6j+k}{2}>$ & $\<7,1+i,\frac{7+7j}{2},\frac{8+i+2j+k}{2}>$ \\
& $\<7,7i,\frac{1+8i+j}{2},\frac{8+i+k}{2}>$ & $\<7,7i,\frac{1+8i+j}{2},\frac{12+5i+k}{2}>$ \\
& $\<7,7i,\frac{3+10i+j}{2},\frac{10+3i+k}{2}>$ & $\<7,7i,\frac{3+10i+j}{2},\frac{8+i+k}{2}>$ \\
$\id{b}_3$
& $\<7,7i,\frac{5+12i+j}{2},\frac{12+5i+k}{2}>$ & $\<7,7i,\frac{5+12i+j}{2},\frac{4+11i+k}{2}>$ \\
& $\<7,7i,\frac{9+2i+j}{2},\frac{2+9i+k}{2}>$ & $\<7,7i,\frac{9+2i+j}{2},\frac{10+3i+k}{2}>$ \\
& $\<7,7i,\frac{11+4i+j}{2},\frac{4+11i+k}{2}>$ & $\<7,7i,\frac{11+4i+j}{2},\frac{6+13i+k}{2}>$ \\
& $\<7,7i,\frac{13+6i+j}{2},\frac{6+13i+k}{2}>$ & $\<7,7i,\frac{13+6i+j}{2},\frac{2+9i+k}{2}>$ \\[1ex]
\hline & & \\[-4ex]
& $\<7,2+i,\frac{7+7j}{2},\frac{9+i+5j+k}{2}>$ & $\<7,2+i,\frac{7+7j}{2},\frac{9+i+j+k}{2}>$ \\
& $\<7,7i,\frac{1+4i+j}{2},\frac{2+i+k}{2}>$ & $\<7,7i,\frac{1+4i+j}{2},\frac{6+3i+k}{2}>$ \\
& $\<7,7i,\frac{3+12i+j}{2},\frac{6+3i+k}{2}>$ & $\<7,7i,\frac{3+12i+j}{2},\frac{4+9i+k}{2}>$ \\
$\id{b}_4$
& $\<7,7i,\frac{5+6i+j}{2},\frac{10+5i+k}{2}>$ & $\<7,7i,\frac{5+6i+j}{2},\frac{2+i+k}{2}>$ \\
& $\<7,7i,\frac{9+8i+j}{2},\frac{4+9i+k}{2}>$ & $\<7,7i,\frac{9+8i+j}{2},\frac{12+13i+k}{2}>$ \\
& $\<7,7i,\frac{11+2i+j}{2},\frac{8+11i+k}{2}>$ & $\<7,7i,\frac{11+2i+j}{2},\frac{10+5i+k}{2}>$ \\
& $\<7,7i,\frac{13+10i+j}{2},\frac{12+13i+k}{2}>$ & $\<7,7i,\frac{13+10i+j}{2},\frac{8+11i+k}{2}>$ \\[1ex]
\hline
\end{tabular}
\caption{Table of left $\Ops$- and $\Ons$-subideals.
          \label{table:O:7}}
\end{table}




\section{Gross formula for level $p^2$}

The aim here is to conjecture a formula that applies to modular forms
of level~$p^2$, similar to the one proved by Gross
in~\cite[Proposition 13.5, p.179]{Gross} for prime level.
Keep the notation of the previous section, and let $f$ be a newform of
weight $2$ and level $p^2$, such that $f$ is an eigenform for the Hecke operators.
We denote by $L(f,s)$ the Hecke $L$-series of $f$,
and for $D$ a fundamental discriminant we define its
\emph{twisted $L$-series} as
\[
    L(f,D,s) := L\left(f\otimes D,s\right),
\]
where $f \otimes D$ is (the newform corresponding to) the
twist of $f$ by the quadratic character $n\mapsto\kro{D}{n}$.
Recall that $L(f,D,s)$ is an entire function of the complex plane with
a functional equation relating its values at $s$ and $2-s$, with
central value $L(f,D,1)$.
The sign of the functional equation will be denoted by
$\epsilon(f,D)$. This sign determines the parity of the order of
vanishing of $L(f,D,s)$ at $s=1$; in particular when
$\epsilon(f,D)=-1$ it is trivial that $L(f,D,1)=0$.

We denote by $\M(\Ot)^f$ the $f$-isotypical component of $\M(\Ot)$.
i.e. the eigenspace of $\HeckeRing^0$
with the same eigenvalues as $f$.
We have the following result due to Pizer (\cite[Theorem 8.2, p.223]{Pizer}) : 
\[
  \dim \M(\Ot)^f = \begin{cases}
     2 & \text{if $f$ is not the twist of a level $p$ form,} \\
     1 & \text{if $f$ is a level $p$ form twisted by a quadratic character,}\\
     0 & \text{if $f$ is a level $p$ form twisted by a
     non-quadratic character.}
  \end{cases}
\]
In what follows we will assume that $f$ is not in the last case, i.e.
that $\M(\Ot)^f \neq 0$.

Since $\level(\Ot)=4p$, it follows that
$\Gross\bigl(\M(\Ot)^f\bigr)=0$. In order to obtain modular forms of
weight $\tf{3}{2}$ and level $4p^2$, we have to employ the orders of
level $p^2$ introduced in \S\ref{sec:ordersp2}.
Let $\Os$ be such an order. The general construction of
\S\ref{sec:quatshimura} gives a Hecke-linear map
\[
   \Gross:\M(\Os)\rightarrow
        M_{\tf{3}{2}}(4p^2,\charp).
\]
where $\charp(n) := \kro{p}{n}$ is the character of $\Os$.
However, the space $\M(\Os)$ represents
weight $2$ modular forms of level $p^3$, and hence it is too
big. Indeed, by Theorem~\ref{thm:idealsp2}~(3) we know that
$\dim\M(\Os)=p\dim\M(\Ot)$.

To overcome this problem, we will use the theory of
\S\ref{sec:ordersp2} and define Hecke-linear maps
\[
   \Gross^{\Ot}_{\Os} := \frac{1}{p} \Theta \circ \psi^{\Ot}_{\Os}
     : \M(\Ot)\rightarrow M_{\tf{3}{2}}(4p^2,\charp),
\]
which in view of Theorem~\ref{thm:idealsp2} can also be defined, for
$\cls{b}\in\I(\Ot)$, by
\[
   \Gross^{\Ot}_{\Os}(\cls{b}) := \Gross(\cls{c}),
\]
where $\id{c}$ is any ideal in $\psi^{\Ot}_{\Os}(\id{b})$.
Note that, although the map $\Gross^{\Ot}_{\Os}$
depends on $\Ot$ and $\Os$, its image depends only on $\sigma(\Os)$,
so that we really have only two different constructions.

\begin{prop}\label{prop:thetazero}
Let $p^\ast = \kro{-1}{p} p$ be the prime discriminant associated to $p$.
Note that the level of $f\otimes p^\ast$ can be $p$ or $p^2$.
If $\epsilon(f,1)=+1$ and either
\begin{itemize}
\item[(A)] $f\otimes p^\ast$ has level $p$
      and $\epsilon(f,p^\ast) = \kro{-1}{p} \sigma(\Os)$, or
\item[(B)] $f\otimes p^\ast$ has level $p^2$
      and $\epsilon(f,p^\ast) = \kro{-1}{p}$,
\end{itemize}
then $\Gross^{\Ot}_{\Os}(\M(\Ot)^f) = 0$.
\end{prop}
\begin{proof}
From $\epsilon(f,1)=+1$ it follows that $\epsilon(f,-d)=-1$ for any
fundamental discriminant $-d<0$ not divisible by $p$.  Either of the
conditions imply that for any fundamental discriminant $-pd<0$ with
$\kro{d}{p}=\sigma(\Os)$, we have $\epsilon(f,-pd)=-1$ (see Lemma 30
and Theorem 6 of \cite{Atkin}) , and in particular $L(f,-pd,1)=0$. The
result now follows from Waldspurger's formula~\cite[Th\'eor\`eme 1]{Waldspurger}.
\end{proof}

It is worth noting that the first examples of $f$ satisfying
condition (B) above
(e.g. the modular form of level $13^2$ and degree $2$ over $\Q$,
as well as the two rational modular forms of level $37^2$ and rank $0$,
or one of the rational modular forms of level $43^2$ and rank $0$),
we found that indeed $\M(\Ot)^f=0$.
It is an interesting question whether condition (B) characterizes all
modular forms of level $p^2$ and trivial character coming from level
$p$ ones by twisting by a non-quadratic character.

\begin{conj} Assume $\M(\Ot)^f\neq 0$
\begin{enumerate}
\item if $L(f,1)=0$, then $\Gross^{\Ot}_{\Os}(\M(\Ot)^f) = 0$.
\item if $L(f,1)\neq 0$, then $\Gross^{\Ot}_{\Os}(\M(\Ot)^f) \neq 0$,
      unless $f$ and $\Os$ are in the conditions (A) or (B) of
      Proposition~\ref{prop:thetazero}.
\end{enumerate}
\end{conj}
\begin{remark}
In case $L(f,1)\neq 0$, the conjecture asserts that our construction
yields two different weight $\tf{3}{2}$ modular forms (in Kohnen space)
in Shimura correspondence with $f$, unless $f\otimes p^\ast$ has level
$p$, where the construction yields only one modular form, or if $f$ is
in condition (B), where we obtain no modular forms.
It follows from~\cite[Prop 3]{Ueda2} that this is exactly what is
available in Kohnen space, even when $L(f,1)=0$.
\end{remark}

Note that the analogue of our conjecture for the case of level $p$ is
true, and follows from Gross's formula~\cite{Gross}. This has been
extended to the case of odd square free level in~\cite{BSP90} and
\cite{BSP92}. To support our conjecture we will give an explicit
version, namely a formula, which has been verified numerically for
many cases.

The strong multiplicity one theorem of Ueda~(\cite[Theorem 3.11,
p.181]{Ueda}) implies that $\dim\Gross^{\Ot}_{\Os}(\M(\Ot)^f)\leq 1$.
The above conjecture gives conditions on $\Os$ and $f$ so that
$\Gross^{\Ot}_{\Os}\bigl(\M(\Ot)^f\bigr)\neq0$.
In that case, $\ker\Gross^{\Ot}_{\Os}$ has codimension $1$
in $\M(\Ot)^f$, and thus
a vector $\e{f,\Os}\in\M(\Ot)^f$ is well
defined up to a constant by requiring it to be orthogonal to
$\ker\Gross^{\Ot}_{\Os}$; we will write
\[
   \Gross_{f,\Os} := \Gross^{\Ot}_{\Os}(\e{f,\Os}) =
         \sum_{d\geq 1} c_{f,\Os}(d)\, q^d.
\]
Otherwise, although our method yields the zero modular form, we will
use $\e{f,\Os}$ to denote any nonzero vector in $\M(\Ot)^f$, with the
understanding that it doesn't really matter which one we choose since
$\Gross_{f,\Os}=0$ regardless of this choice, because our conjectural
formula below is nontrivial even in this case.

Define the rational constant $\alpha_f$ by
\[
  \alpha_f := \frac{1}{2}\cdot\begin{cases}
     1             & \text{if $f$ is not the twist of a level $p$ form,} \\
     \frac{p}{p-1} & \text{if $f$ is the quadratic twist of a level $p$ form.}\\
  \end{cases}
\]

\begin{conj}
\label{conj:formula}
Let $d$ be an integer such that $-pd<0$ is a
fundamental discriminant, and such that $\kro{d}{p}=\sigma(\Os)$. Then
\[
  L(f,-pd,1)\, L(f,1)
  = \alpha_f\,\frac{\<f,f>}{\sqrt{pd}}
     \frac{c_{f,\Os}(d)^2}{\<\e{f,\Os},\e{f,\Os}>}.
\]
\end{conj}
\begin{remark}
The formula is true up to a constant (depending on $f$ and $\Os$), by
Waldspurger's theorem.  Moreover, in case $f$ is rational it is possible to
prove that
$L(f,-pd,1)L(f,1)\sqrt{pd}/\<f,f>$ is a rational number of bounded
height, and a similar result holds for general $f$.
Thus there is an effective procedure to determine if the
formula is true for any particular $f$.
See \cite{Pa-To} for some numerical examples in the cases
$p=7,11,17,19$, where the formula has been verified.
\end{remark}

We note the following complement of Proposition~\ref{prop:thetazero}
\begin{lemma}
Assume $f$ and $\Os$ are such that neither of the conditions (A) nor
(B) of Proposition~\ref{prop:thetazero} hold. Then there is a
fundamental discriminant $-pd<0$ such that $\kro{d}{p}=\sigma(\Os)$,
and $L(f,-pd,1)\neq 0$.
\end{lemma}
\begin{proof}
It is easy to see that if neither (A) nor (B) hold, then for any
fundamental discriminant $-pd'<0$ such that $\kro{d'}{p}=\sigma(\Os)$,
we have $\epsilon(f,-pd')=+1$. It now follows from~\cite{Bump} that there
is a fundamental discriminant $d''$ prime to $d'$ such that $d''>0$
and $\kro{d''}{p}=+1$, for which $L(f\otimes-pd', d'',1)\neq 0$.
It is clear that $d=d'd''$ will satisfy the claimed conditions.
\end{proof}

\begin{prop}
Conjecture 2 implies Conjecture 1.
\end{prop}
\begin{proof}
If $L(f,1)=0$, then clearly the formula implies that $\Gross_{f,\Os}=0$.
Conversely, if the conditions of
Proposition~\ref{prop:thetazero} are not satisfied, the Lemma asserts
the existence of some $d$ in the hypothesis of Conjecture 2 for which
$L(f,-pd,1)\neq 0$. Thus, if $L(f,1)\neq 0$, it follows that
$c_{f,\Os}(d)\neq 0$, i.e. $\Gross_{f,\Os}\neq 0$.
\end{proof}


\begin{remark}
  When $p\equiv 3\pmod{4}$, and $f$ is a modular form of level $p$ or
  $p^2$, this method can be used to compute a weight $\tf{3}{2}$
  modular form whose Fourier coefficients are related to the central
  values of \emph{real} quadratic twists of $f$.
  See \cite{Pa-To} for an exposition of this method and examples.
\end{remark}



\begin{thebibliography}{B\"o-SP90}

\bibitem[At-Le]{Atkin} A. Atkin and J. Lehner,
  \emph{Hecke Operators on $\Gamma_0(m)$},
  Math. Ann. 185 (1970), p. 134-160

\bibitem[B\"o-SP90]{BSP90} S. B\"ocherer and R. Schulze-Pillot,
  \emph{On a theorem of Waldspurger and on Eisenstein series of
  Klingen type},
  Math. Ann. 288 (1990) 361--388.

\bibitem[B\"o-SP92]{BSP92} S. B\"ocherer and R. Schulze-Pillot,
  \emph{The Dirichlet series of Koecher and Maass and modular forms
         of weight $\frac32$},
  Math. Z.  209 (1992),  no. 2, 273--287.

\bibitem[BFH]{Bump} D. Bump, S. Friedberg and J. Hoffstein
  \emph{Nonvanishing theorems for L-functions of modular forms and
  their derivatives},
Invent. Math. V. 102, 1990, p. 543-618.

\bibitem[Gr]{Gross} B. Gross, \emph{Heights and Special Values of
  L-series},
Conference Proceedings of the CMS Vol. 7, 115-187
AMS 1986

\bibitem[Ll]{Llorente} P. Llorente,
\emph{Correspondencia entre formas ternarias enteras y \'ordenes
cuaterni\'onicos},
Rev. R. Acad. Cienc. Exactas F\'{\i}s. Nat. (Esp.)
{\bf 94}, no. 3, 2000, 397--416.

\bibitem[Pa-Vi]{Ville}  A. Pacetti and F. Rodriguez-Villegas, \emph{Computing
  Weight $2$ Modular forms of level $p^2$},
  Math. Comp. {\bf 74} (2005), 1545--1557.

\bibitem[Pa-To]{Pa-To} A. Pacetti and G. Tornar\'\i a, \emph{Examples of
  Shimura correspondence for level $p^2$ and real quadratic twists},
to appear in ``Ranks of elliptic curves and random matrix theory.''
Available online at arxiv.org/math.NT/0412104

\bibitem[Pi]{Pizer} A. Pizer, \emph{Theta Series and Modular Forms of Level
$p^2M$}
Compositio Mathematica, Vol. {\bf 40} Fasc. 2, 1980, p. 177--241.

\bibitem[Sh]{Shimura} G. Shimura \emph{On Modular Forms of Half
Integral Weight}
The Annals of Mathematics, Vol. 97, 1971, p. 440-481.


\bibitem[To]{ell-data}
Gonzalo Tornar\'{\i}a, \emph{Data about the central values of the {L}-series of
  (imaginary and real) quadratic twists of elliptic curves}, {\tt
  http://www.ma.utexas.edu/users/tornaria/cnt/}, 2004.



\bibitem[Ue88]{Ueda2} M. Ueda,
  \emph{The decomposition of the spaces of cusp forms of half-integral
  weight and trace formula of Hecke operators},
  J. Math. Kyoto Univ. 28-3 (1988) 505--555.


\bibitem[Ue93]{Ueda} M. Ueda,
  \emph{On twisting operators and newforms of half-integral weight},
Nagoya Mathematical Journal Vol. 131,1993, p. 135-205.


\bibitem[Wa]{Waldspurger} J-L Waldspurger \emph {Sur les coefficients de Fourier des
  formes modulaires de poids demi-entier}
Journal de Mathematiques Pures et Appliquees, Tome 60 Fas. 4, 1981 p. 375-484

\end{thebibliography}
\end{document}